\documentclass[11pt,a4paper]{article}
\usepackage{graphicx,amssymb,amsfonts,latexsym,amsmath,amsthm,times}
\usepackage{epsfig}
\usepackage{fancyhdr}
\usepackage{color}
\usepackage[english]{babel}
\usepackage{authblk}
\usepackage{fancyhdr}
\usepackage{url}
\usepackage{tabularx}
\usepackage{float}
\usepackage{hyperref}
\usepackage{rotating}
\usepackage{graphics}
\usepackage{hyperref}
\usepackage{extarrows}
\usepackage{xcolor}
\usepackage{lscape}
\usepackage{multirow}
\usepackage{afterpage}
\usepackage{longtable}
\usepackage[T1]{fontenc}
\usepackage{calligra}
\usepackage{url}  

\newcommand\overcirc[1]{\raisebox{10pt}{\tiny$\circ$}{\kern-6.5pt}\mbox{$#1$}}
\newcommand\undersym[2]{\raisebox{-5pt}{\tiny$#2$}{\kern-12pt}\mbox{$#1$}}
\newcommand\oversym[2]{\raisebox{8pt}{\tiny$#2$}{\kern-10pt}\mbox{$#1$}}
\newtheorem{thm}{Theorem}[section]

\begin{document}
\title{Exact Solutions of the Time Derivative Fokker-Planck Equation: A Novel Approach }
\author{H. I. Abdel-Gawad$^{1}$$^{*}$, \ \  N. H. Sweilam$^{1}$, \ \ S. M. AL$-$Mekhlafi$^{2}$, \ \ D. Baleanu$^{3,4}$ }
\date{
$^{1}$Mathematics Department, Faculty of Science Cairo University, Egypt.\\
$^{2}$Mathematics Department, Faculty of Education, Sana'a University, Yemen\\
	$^{3}$  Cankaya University, Department of Mathematics, Turkey\\
$^{4}$Institute of Space Sciences, Magurele-Bucharest, Romania\\
 E-mail:  hamdyig@yahoo.com$^{1,\ast}$, \ \ \ nsweilam@sci.cu.edu.eg $^{1}$,\,\,\ smdk100@yahoo.com$^{2}$,\,\,\ 
 dumitru@cankaya.edu.tr$^{3,4}$
}

\maketitle
\textbf{Abstract:}
In the present article, an approach to find the exact solution of the fractional
Fokker-Planck equation is presented. It is based on transforming
it to a system of first-order partial differential equation via Hopf transformation, together
with implementing the extended unified method. On the other hand,
reduction of the fractional derivatives to non autonomous
ordinary derivative. Thus the fractional Fokker-Planck equation is reduced
to non autonomous classical ones. Some explicit solutions of the classical,
 fractional time derivative Fokker-Planck equation, are obtained . It is shown
that the solution of the Fokker-Planck equation is bi-Gaussian's. It is found that high
friction coefficient plays a significant role in lowering the standard
deviation. Further, it is found the fractionality has stronger effect
than fractality. It is worthy to mention that the mixture of Gaussian's
is a powerful tool in machine learning. Further, when varying the
order of the fractional time derivatives, results to slight
effects in the probability distribution function. Also, it is shown
that  the mean and mean square of the velocity vary slowly.
\begin{flushleft}
\end{flushleft}
\section{Introduction}
The Fokker-Planck equation  (FPE) deals with fluctuations of systems which stem from disturbances,
each of which changes the variables of the system in an unpredictable
way. When macroscopic particles are immersed in fluid, they are pulled
by the fluid and the position of particles is an unpredictable. So
that they fluctuate about an expected position with a certain probability
to find the particle in a given region. The probability density function
can be determined via the FPE.This equation is used in many different
fields in natural science, in solid-state physics, quantum optics,
chemical physics, and mathematical biology. Anomalous diffusion is
an eminent behavior and characterizes the transport processes of problem
in several models {[}1\textendash 4{]}.The features of diffusion
can be classified to sub- and super- diffusion. They may be measured
by evaluating the variance, relevant to the displacement of particles
in the medium, which behaves as $x^{\beta}$where $x$ is the space variable.
When $0<\beta<1$, this stands to sub-advection (or in transition
state) and when $1<\beta<2$, this case is sub-diffusion (or the transition from advection
to diffusion). While in the case when $2<\beta<3$ describes super-diffusion
(or the transition from diffusion to dispersion). The first case assigns
to spatially disordered or fractal media and also in fractional Brownian
motion and random processes [5, \ 6]. To this end the fractal and
fractional time-space derivatives, of orders $\alpha,\,\,0<\alpha<1$
and $\beta,\:1<\beta<2,$ respectively, do emerge in the formulation
of diffusion equations relevant to the study of the transport process.
By using the time time fractional derivative, the effects of time
distributed delay in the transport process are then taken into consideration.
While fractal derivative reflect instantaneous frangibility effect.

Fractional derivatives which was proposed in the literature assume
numerous definitions; Caputo, Caputo-Fabrizio [7, \ 8], and a very
recently Atangana-Baleneau derivatives {[}9{]} and slso Riemann-Liouville,
Riesz (see {[}21{]}), In the last three, only the kernel of the integral
is different depending to when it is singular or not. Analytical investigation
of the fractional KPE equation that describes anomalous diffusion
of energetic particles {[}10{]} are accomplished. An approach
was proposed in {[}11\textendash 13{]} to solve analytically the fractional
force-less 1-D FPE. The Klein\textendash Kramers equation (KKE)
{[}14,15{]} that describes the transport of energetic particles in
turbulent magnetic fields are often reduced  to the force-less homogeneous
one-dimensional FPE.

The KKE {[}33, \ 34{]}, which depicts the 
the Brownian motion of particles in the presence of an outer force
$F(z)$ is

\begin{equation}
\frac{\partial}{\partial t}W(t,z,v)=(-v\frac{\partial}{\partial z}+\frac{\partial}{\partial v}(\eta v-\frac{F(z)}{m})+B\frac{\partial^{2}}{\partial v^{2}})W(t,z,v),
\end{equation}
where $W(t, z,v)$ is the probability density function of particles,
$B=\frac{\eta K_{B}T}{m}$, \ \ $K_{B}$ is the Boltzmann constant, $T$
is the absolute temperature and $m$ is the mass of a diffusing particle,
$v$ its velocity, $\eta$ denotes the friction coefficient, $F(z)$
is an external force field.

To depict the distribution in the velocity space, we evaluate
the mean and mean square of the particle position in the absence of
the external force. In this case, (1) reduces to the FPE that inspect
the diffusion of a test particle in the phase space,

\begin{equation}
\frac{\partial}{\partial t}f(t,v)=\eta\frac{\partial}{\partial v}(f(t,v)v)+B\frac{\partial^{2}}{\partial v^{2}}f(t,v),
\end{equation}
 $f(t,v)=\int_{-\infty}^{\infty}W(t,z,v)dz$ designates the distribution
function.

Here, we consider the  fractional time KPE of (2) that
generalizes the energetic particles transport equations, that is by
replacing the ordinary derivatives using the  Caputo and
Caputo-Fabrizio fractional time derivatives {[}7,8{]} to create a prediction
of the evolution of the particle distribution function in phase space.
We mention that the results obtained incorporate the effects of the
distributed time delay.

%
%
%
%
%
%
%

\section{Fractional Derivative}

The Caputo fractional derivative (CFD) is defined as follows [21]:

\begin{equation}
_{0}^{C}D^{\alpha}_{t}f(t)=\frac{1}{\Gamma(m-\alpha)}\int_{0}^{t}(t-t_{1})^{-\alpha+(m-1)}f^{(m)}(t_{1})dt_{1},\:m-1<\alpha<m, \ \  \ t>0,
\end{equation}
provided that $f$ is Hoder continuous , $f\epsilon H^{m,\alpha-(m-1)}(\mathbb{R}^{+})$
and the integral exists.\\

 The Caputo-Fabrizio fractional derivative
(CFFD)  is  defined as follows [7]:

\begin{equation}
_{0}^{CF}D^{\alpha}_{t}f(t)=\frac{2\alpha}{(1-\alpha)(2-\alpha)}\int_{0}^{t}e^{-\frac{\alpha}{1-\alpha}(t-t_{1})}f^{\prime}(t_{1})dt_{1},\ \ 0<\alpha<1, \ \ \ t>0,
\end{equation}
provided that $f\epsilon H^{1,\alpha}(\mathbb{R}^{+})$.

The Atangana-Baleanu fractional derivative (ABFD), in the Caputo sense,
is  defined as follows {[}9,20{]}:

\begin{equation}
_{0}^{AB}D^{\alpha}_{t}f(t)=\frac{B(\alpha)}{1-\alpha}\int_{0}^{t}E_{\alpha}(-\frac{\alpha}{1-\alpha}(t-t_{1})^{\alpha})f^{\prime}(t_{1})dt_{1},\ \ \ 0<\alpha<1, \ \ \ t>0,
\end{equation}
where $B(\alpha)>0$ is a normalization function satisfying $B(0)=B(1)=1$
 and $E_{\alpha}(t)$ is the Mittag\textendash Leffler function and
$f\epsilon H^{1,\alpha}(\mathbb{R}^{+})$. It is worth noticing that
this function is not invariant under the CFD. The function which is
invariant is $e_{\alpha}(t)$ {[}26{]} where

\begin{equation}
E_{\alpha}(t)=\sum_{n=0}^{\infty}\frac{t^{n}}{\Gamma(\alpha n+1)}, \ \ \ e_{\alpha}(t):=e_{\alpha,1}(t)=\sum_{n=0}^{\infty}\frac{t^{\alpha n}}{\Gamma(\alpha n+1)},
\end{equation}
and the last fractional exponential function generalizes to

\begin{equation}
e_{\alpha,\beta}(\lambda,t)=\sum_{n=0}^{\infty}\frac{\lambda^{n}t^{\alpha n}}{\Gamma(\alpha n+\beta)}.
\end{equation}

Here, we propose a new fractional derivative called  Gawad's fractional derivative  (GFD) which is defined  as follows:

\begin{equation}
_{0}^{G}D^{\beta}_{t}f(t)=\frac{2\lambda^{\frac{1}{\beta}}}{\lambda+2}\int_{0}^{t}e^{-\lambda(t-t_{1})^{\beta}}f^{\prime}(t_{1})dt_{1},\  \ \beta>0, \ \ \ t>0,
\end{equation}
provided that $f\epsilon H^{1,\beta}(\mathbb{R}^{+})$. We mention
that the CFFD is a particular case from (8), when $\beta=1$and $\lambda=\frac{\alpha}{1-\alpha}$. The
kernel in the fraction derivative (8) is of interest) in the theory
of distributions. As when $0<\beta\leq1$, the fractional exponential
distribution is $f(t)=\lambda^{1/\beta}e^{-\lambda t^{\beta}}$generalizes
the classical exponential distribution and when $1<\beta\leq2$, the
fractional Gaussian distribution is $f(t)=\frac{\lambda^{1/\beta}}{\sqrt{\pi}}e^{-\lambda t^{\beta}}.$
This later distribution may be considered as transition between the
exponential and the Gaussian's.

\subsection{Reduction of the Fractional Derivatives}

Here, we shall reduce the fractional derivatives to non autonomous ordinary derivatives,
that is with time dependent coefficients {[} 22-24{]}.

We consider (3) when $0<\alpha<1,$  is given as follows:

\begin{equation}
_{0}^{C}D^{\alpha}_{t_{1}}f(t_{1})=\frac{1}{\Gamma(1-\alpha)}\int_{0}^{t_{1}}(t_{1}-t_{2})^{-\alpha}f^{\prime}(t_{2})dt_{2}.
\end{equation}
\begin{thm}
 We assume that the integrand in (9 ) is uniformly integrable
on $[0,T_{0}],$  then the CFD is reduced to

\begin{equation}
_{0}^{C}D^{\alpha}_{t}f(t)=\frac{1}{\Gamma(2-\alpha)}(T_{0}-t)^{1-\alpha}f^{\prime}(t),\ \ \ 0\leq t\leq T_{0}.
\end{equation}
\end{thm}
$Proof$ In (9), when operating by the integral, over $t_{1}$ on
$[0,t]$, it holds that

\begin{equation}
\int_{0}^{t}{_{0}^{C}D^{\alpha}_{t_{1}}f(t_{1})}\,dt_{1}=\int_{0}^{t}\Bigg(\frac{1}{\Gamma(1-\alpha)}\int_{0}^{t_{1}}(t_{1}-t_{2})^{-\alpha}f^{\prime}(t_{2})dt_{2}\Bigg)dt_{1}.\ \ 
\end{equation}
By using the assumption, we permute the inner with the outer integral
in the RHS of (11), and have

\[
\frac{1}{\Gamma(1-\alpha)}\int_{0}^{t}\Bigg(\int_{t_{2}}^{t}(t_{1}-t_{2})^{-\alpha}dt_{1}\Bigg)f^{\prime}(t_{2})dt_{2},
\]
or

\begin{equation}
\int_{0}^{t} {_{0}^{C}D^{\alpha}_{t_{1}}} f(t_{1})dt_{1}=\frac{1}{\Gamma(1-\alpha)}\int_{0}^{t}\frac{(t-t_{2})^{1-\alpha}}{1-\alpha}f^{\prime}(t_{2})\,dt_{2}.
\end{equation}
From (12), it holds that

\begin{equation}
_{0}^{C}D^{\alpha}_{t}f(t)=\frac{1}{\Gamma(2-\alpha)}(T_{0}-t)^{1-\alpha}f^{\prime}(t), \ \ \ 0\leq t\leq T_{0}.
\end{equation}
Which is identically (9).$\square$

By using the transformation $f(t):=\tilde{f}(\tau),$ (16) can be
rewritten:

\begin{equation}
_{0}^{C}D^{\alpha}_{t}f(t):=\frac{d}{d\tau}\tilde{f}(\tau),\ \ \tau=\frac{\Gamma(2-\alpha)}{\alpha}(T_{0}^{\alpha}-(T_{0}-t)^{\alpha}).
\end{equation}

The equations (10) and (14) are among the main results in this work.

By the same way, we find

\begin{equation}
\begin{array}{c}
_{0}^{CF}D^{\alpha}_{t}f(t)=\frac{2}{(2-\alpha)}\Bigg(1-e^{-\frac{\alpha}{1-\alpha}(T_{0}-t)}\Bigg)f^{\prime}(t):=\frac{d}{d\tau}\tilde{f}(\tau),\\
\tau=\frac{(2-\alpha)(1-\alpha)}{2\alpha}Log\Bigg(\frac{e^{\frac{\alpha}{1-\alpha}(T_{0}-t)}-1}{e^{\frac{\alpha}{1-\alpha}T_{0}}-1}\Bigg),\ \ 0\leq t\leq T_{0},
\end{array}
\end{equation}

\begin{equation}
_{0}^{AB}D^{\alpha}_{t}f(t)=\frac{B(\alpha)}{(1-\alpha)}(T_{0}-t)e_{\alpha,2}\Bigg(-\frac{\alpha}{1-\alpha},(T-t_{1})\Bigg)f_{0}^{\prime}(t),
\end{equation}
where {[}25{]}

\begin{equation}
e_{\alpha,\beta}(\sigma,x)=\sum_{n=0}^{\infty}\frac{\sigma^{n}x^{\alpha n}}{\Gamma(\alpha n+\beta)},
\end{equation}
and

\begin{equation}
\begin{array}{c}
_{0}^{G}D^{\beta}_{t}f(t)=\frac{2}{\lambda+2}\gamma(\frac{1}{\beta},\lambda(T_{0}-t)^{\beta})f^{\prime}(t):=\frac{d}{d\tau}\tilde{f}(\tau),\\
\tau=\frac{(\lambda+2)}{2}\int_{0}^{t}\frac{1}{\gamma(\frac{1}{\beta},\lambda(T_{0}-t_{1})^{\beta})}dt_{1},\ \ 0\leq t\leq T_{0}, \ \ 0<\beta<1,
\end{array}
\end{equation}
where $\gamma(m,t)$ is the incomplete lower Gamma function. We mention
that

\begin{equation}
\gamma(a,x)+\varGamma(a,x)=\Gamma(a),\ \ \ a>0,\ \ \ x>0,\ \ \ \ \varGamma(a,x)=\int_{x}^{\infty}e^{-y}y^{a-1}dy.
\end{equation}
Now, we will prove the following theorem by using (19).
\begin{thm}
The GFD satisfies the following:

(i) $_{0}^{G}D^{\beta}_{t}(g(t)+f(t))=_{0}^{G}D^{\beta}_{t}g(t)+_{0}^{G}D^{\beta}_{t}f(t)$.\\

(ii) $_{0}^{G}D^{\beta}_{t}(f(t)g(t))=f(t)$$_{0}^{G}D^{\beta}_{t}g(t)+g(t)_{0}^{G}D^{\beta}_{t}f(t).$\\

(iii)$_{0}^{G}D^{\beta}_{t}(\frac{f(t)}{g(t)}))=\frac{g(t)_{0}^{G}D^{\beta}_{t}f(t)-f(t)_{0}^{G}D^{\beta}_{t}g(t)}{g(t)^{2}}$.

Further, by using (19) the function $f(t)$ which is invariant under
the GFD and   
$$_{0}^{G}D^{\beta}_{t}e_{\beta,G}(t)=e_{\beta,G}(t),$$
is found directly by:
\begin{equation}
e_{\beta,G}(t)=e^{\frac{(\lambda+2)}{2}\int_{0}^{t}\frac{1}{\gamma(\frac{1}{\beta},\lambda(T_{0}-t_{1})^{\beta})}dt_{1}}.
\end{equation}
\end{thm}

\section{Solutions of the FPE}

We present a novel approach to solve of (2) as follows. Here we use
the transformations $$f_{v}(t,v)=F(v,t)\,f(t,v), \ \ \  \  f_{t}(t,v)=G(t,v)\,f(t,v).$$
Thus (2) is rewritten as follows:
\begin{equation}
\begin{array}{c}
f_{v}(t,v)-F(t,v)\,f(t,v)=0,\ \ \ f_{t}(t,v)-G(t,v)\,f(t,v)=0,\\
G(t,v)-\eta(1+vF(t,v))-B(F_{v}(t,v)+F(t,v)^{2})=0.
\end{array}
\end{equation}

On the other hand, the unified method is implemented. It asserts
that the solutions of  partial differential equations (PDEs) can be expressed by a polynomial or a rational
forms in an auxiliary functions that satisfy appropriate auxiliary
equations. Here, we are concerned with rational forms in the cases
of different auxiliary equations. This will be classified in what
follows.

\subsection{Case of Linear Auxiliary Equations}

We write the solutions of (21) in the form.

\begin{equation}
\begin{array}{c}
f(t,v)=\frac{s_{1}(v)g(t,v)+s_{0}(v)}{a_{1}(v)g(t,v)+a_{0}(v)},\ \ \ F(t,v)=\frac{b_{1}(v)g(t,v)+b_{0}(v)}{s_{1}(v)g(t,v,z)+s_{0}(v)},\\
G(t,v)=\frac{d_{1}(v)g(t,v)+d_{0}(v)}{s_{1}(v)g(t,v)+s_{0}(v)}.
\end{array}
\end{equation}

We assume that the auxiliary equations are linear, 

\begin{equation}
\begin{array}{c}
g_{t}(t,v)=\mu\,(c_{1}g(t,v)+c_{0}),\ \ \ \ g_{v}(t,v)=h(v)\,(c_{1}g(t,v)+c_{0}).\end{array}
\end{equation}

It is worth noticing that in (23), the compatibly equation $g_{tv}(t,v)=g_{vt}(t,v)$
holds.

By inserting (22) and (23) into (21), and by using (24), we get a
system of nonlinear algebraic equations in \ $a_{i}(v), \ \ b_{i}(v) \ and \ d_{i}
(v),$\ \ \ $i=1,2... \ .$ \ We found that the calculations are not straight forward
due to two facts.\\
(i)The equations obtained are nonlinear.\\
(ii) It arises that we can have two equations in the same dependent
variable, for example, \ $a_{j}(v)$ \ and  \ $a_{j}^{\prime}(v),\ \ \ j=0,1,... ,$
so, we have to use the compatibly equation, $b_{j}^{\prime}(v)-(b_{j}(v))^{\prime}=0$.

By setting the coefficients of $g(t,v)^{j}, \ \ \ j=0,1,... ,$ equal to zero,
we have the following equations
\begin{align}
a_{0}^{\prime}(v)&=\frac{1}{s_{0}(v)}\Bigg(-c_{0}a_{1}(v)h(v)s(v)+a_{0}(v)\Bigg(-b_{0}(v)\notag\\
&+c_{0}h(v)s_{1}(v)+s_{0}^{\prime}(v)\Bigg)\Bigg),\notag\\
s_{1}^{\prime}(v)&=b_{1}(v), \notag\\
s_{0}^{\prime}(v)&=b_{0}(v)+h(v)\Bigg(-c_{0}s_{1}(v)+c_{1}s_{0}(v)\Bigg),
\end{align}
\begin{align}
b_{0}^{\prime}(v)&=\frac{1}{m\,B}\Bigg(b_{0}(v)(-m\,v\,\eta+F_{0}+Bc_{1}mh(v))\notag\\
&-m \Bigg(Bc_{0}b_{1}(v)h(v)-c_{0}(k_{0}v+\mu)s_{1}(v)\notag\\
&+\Bigg(c_{1}k_{0}v+\eta+c_{1}\mu \Bigg)s_{0}(v)\Bigg)\Bigg),\notag\\
b_{1}^{\prime}(v)&=\frac{1}{m\,B}(b_{1}(v)(-m\,v\,\eta+F_{0})-m\,\eta s_{1}(v)),\notag
\end{align}
and
\begin{equation}
\begin{array}{c}
d_{1}(v)=0,\;d_{0}(v)=\frac{c_{1}\mu}{a_{1}(v)}\Bigg(s_{1}(v)-a_{1}(v)s_{0}(v)\Bigg),\\
\;a_{1}(v)=c_{1}/c_{0},\ \ \  a_{0}(v)=1.
\end{array}
\end{equation}

By rewriting the second equation in (24) as \  $b_{1}(v)=s_{1}^{\prime}(v),$ \ the
compatibly equation  $$b_{1}^{\prime}(v)-(b_{1}(v))^{\prime}=0,$$ gives
rise to

\begin{equation}
m\,\eta s_{1}(v)+mv\eta s_{1}^{\prime}(v)+B\,ms_{1}^{\prime\prime}(v)=0,
\end{equation}

where (26) solves to

\begin{equation}
\begin{array}{c}
s_{1}(v)=A_{2}e^{-\frac{v^{2}\eta}{2\,B}}+A_{1}\sqrt{\frac{\pi\,B}{2\eta}}e^{\frac{v^{2}\eta}{2\,B}-\frac{F_{0}^{2}}{2Bm^{2}\eta}}erfi(\frac{m\,v\,\eta}{\sqrt{2B\eta}m}),\\
erfi(x)=e^{-x^{2}}\int_{0}^{x}e^{y^{2}}dy.
\end{array}\;
\end{equation}

Thus, we can evaluate $b_{1}(v)=s_{1}^{\prime}(v).$ \\ By the same way,
we rewrite the third equation in (24) as

\begin{equation}
b_{0}(v)=-s_{0}^{\prime}(v)+h(v)(-c_{0}s_{1}(v)+c_{1}s_{0}(v)),
\end{equation}

we find that the compatibly equation $b_{0}^{\prime}(v)-(b_{0}(v))^{\prime}=0,$
gives rise to:
\begin{equation}
(m\,v\,\eta)h(v)+Bc_{1}mh(v)^{2}+m(+\mu-Bh^{\prime}(v))=0.
\end{equation}

The solutions of (29) gives rise to :
\begin{equation}
\begin{array}{c}
c_{1}=\frac{n\,\eta}{\mu},\  \  \ h(v)=\frac{P_{1}(v)}{Q_{1}},\ \ \ P_{1}(v)=-\sqrt{\frac{2B}{\eta}}\mu B_{0}H_{n-1}\Bigg(\frac{v\sqrt{\eta}}{\sqrt{2B}}\Bigg)\\
+v\mu_{1}F_{1}\Bigg(1-\frac{n}{2},\frac{3}{2},\frac{v^{2}\eta}{2B}\Bigg),\\
Q_{1}(v)=BB_{0}H_{n}(\frac{v\sqrt{\eta}}{\sqrt{2B}})+{}_{1}F_{1}\Bigg(-\frac{n}{2},\frac{1}{2},\frac{\eta v^{2}}{2B}\Bigg),
\end{array}
\end{equation}

where $H_{n}(v\sqrt{\frac{\eta}{2B}})$ and $_{1}F_{1}(-\frac{n}{2},\frac{1}{2},\frac{\eta v^{2}}{2B})$
are the Hermite polynomial and the hyper geometric functions respectively.

Further the solution of the auxiliary equations give rise to

\begin{equation}
g(v,t)=-\frac{c_{0}}{c_{1}}+B_{3}e^{c_{1}\int h(v)dv+\mu t}.
\end{equation}

We mention that$\int h(v)dv$ can not be directly evaluated. To this
end we assume that (cf. (30))

\begin{equation}
(Q_{1}+s(v))^{\prime}=P_{1},
\end{equation}

and the calculations give

\begin{equation}
\begin{array}{c}
s(v)=\frac{(n\eta+\mu)}{n\eta}(-BB_{0}H_{n}\Bigg(\frac{v\sqrt{\eta}}{\sqrt{2B}}\Bigg)+B(1-_{P}F_{Q}(\{-\frac{n}{2}\},\{\frac{1}{2}\},\frac{\eta v^{2}}{2B})),\\
\end{array}
\end{equation}
where $_{P}F_{Q}$ is the generalized hyper geometric function. Finally,
we get

\begin{equation}
\int h(v)dv=Log(\mid P_{1}(v)+s(v)\mid).
\end{equation}

By substituting for $s_{i}(v),a_{i}(v)$ and by using (24)-(34) in
the first equation in (22) we get the required solution.

$$f(t,v)=\frac{P}{Q},$$
$$P=c_{0}\mu \Bigg (e^{-\eta\frac{(1+m^{3}v^{2})}{2Bm^{2}}}(2B_{1}e^{\frac{\eta}{2Bm^{2}}})\eta^{3/2}-\sqrt{\frac{2B}{\pi}}\Bigg(2A_{1}Bc_{0}e^{\frac{\eta}{2Bm^{2}}}\Bigg)\Bigg)m-B_{2}\eta)erfi(\frac{m\,v\,\eta}{\sqrt{2B\eta}m}))$$
$$+\eta^{3/2}e^{-\frac{v^{2}\eta}{2B}}(2A_{2}+A_{1}\sqrt{\frac{2B}{\pi\eta}}erfi(\frac{m\,v\,\eta}{\sqrt{2B\eta}m})g(t,v)),$$
\begin{equation}
Q=2\eta^{3/2}(c_{0}\mu+n\eta g(t,v)),
\end{equation}
and
\begin{align}
g(t,v)&=-\frac{c_{0}\mu}{n\eta}+A_{0}e^{n\,t\,\eta+Log(K)}\notag\\
K&=\frac{1}{n\eta}B\Bigg(n\,\eta+\mu-B_{0}\mu H_{n}(v\sqrt{\frac{\eta}{2B}})+n\,\eta\,_{1}F_{1}(-(\frac{n}{2}),\frac{1}{2},\frac{\eta v^{2}}{2B})\notag\\&-(n\,\eta+\mu)_{P}F_{Q} \Bigg(\{-(\frac{n}{2})\},\{\frac{1}{2}\},\frac{\eta v^{2}}{2B}\Bigg)\Bigg).
\end{align}

\subsection{Case of Quadratic Auxiliary Equations}

We assume that the auxiliary equations are quadratic. We write the
solution of (21) in the form

\begin{equation}
\begin{array}{c}
f(t,v)=\frac{s_{1}g(t,v)+s_{0}}{a_{1}(v)g(t,v)+a_{0}(v)},\:F(t,v)=\frac{b_{1}(v)g(t,v)+b_{0}(v)}{s_{1}g(t,v)+s_{0}},\:\\
G(t,v)=\frac{d_{1}(v)g(t,v)+d_{0}(v)}{s_{1}g(t,v)+s_{0}},
\end{array}
\end{equation}

together with auxiliary equation

\begin{equation}
g_{t}(t,v)=\mu\,(c_{2}g(t,v)^{2}+c_{1}g(t,v)+c_{0}),\:g_{v}(t,v)=h(v)\,(c_{2}g(t,v)^{2}+c_{1}g(t,v)+c_{0}).
\end{equation}

It is worth noticing that in (38), the compatibility equation $g_{tv}(v,t)=g_{vt}(v,t)$
holds.

By inserting (35) into (21) and by using (38), again we find that
the calculations are not straight forward due to two facts.

(i) The equations obtained are nonlinear

(ii) It arises that we can have two equations for $a_{j}(v)$ \ and \
$a_{j}^{\prime}(v),\ \ j=0,1,$ \  so we have to use the compatibly equation,
$a_{j}^{\prime}(v)-(a_{j}(v))^{\prime}=0$.
By the same way as in subsection 4.1, we have the following equations
\begin{align}
a_{0}^{\prime}(v)&=\frac{1}{s_{0}}\Bigg(-ao(v)bo(v)-c_{0}s_{0}a_{1}(v)h(v)+c_{0}s_{1}a_{0}(v)h(v)\Bigg),\notag\\
a_{1}^{\prime}(v)&=\frac{1}{s_{1}}\Bigg(-a_{1}(v)b_{1}(v)-c_{2}s_{0}a_{1}(v)h(v)+c_{2}s_{1}a_{0}(v)h(v)\Bigg),\notag\\
a_{1}(v)&=\frac{a_{0}(v)\Bigg((c_{1}^{2}-k_{0}^{2})s_{1}\mu-4c_{2}d_{0}(v)\Bigg)}{\Bigg(c_{1}^{2}-k_{0}^{2}\Bigg)s_{0}\mu},\notag
\end{align} 
\begin{align}
d_{1}(v)&=-\frac{\Bigg(d_{0}(v)\Bigg((c_{1}^{2}-k_{0}^{2})s_{1}\mu-4c_{1}c_{2}s_{0}\mu-4c_{2}d_{0}[v]\Bigg) \Bigg)}{\Bigg(c_{1}^{2}-k_{0}^{2}\Bigg)s_{0}\mu},\notag\\
b_{1}(v)&=\frac{1}{4c_{2}s_{0}}(4c_{2}s_{1}b_{0}(v)-(c_{1}^{2}-k_{0}^{2})s{}_{1}^{2}h(v)\notag\\
&+4c_{1}c_{2}s_{1}s_{0}h(v)-4c_{2}{}_{2}^{2}so^{2}h(v),\notag
\end{align}
\begin{align}
d_{0}(v)&=-\frac{1}{4c_{2}}\Bigg((c_{1}+k_{0})(-c_{1}s_{1}+k_{0}s_{1}+2c_{2}s_{0})\mu\Bigg),\notag\\
b_{0}^{\prime}(v)&=-\frac{1}{4B^{2}c_{2}h(v)^{2}}(B^{2}c_{2}k_{0}^{2}s_{0}h(v)+c2s_{0}(\mu-Bh^{\prime2}(v)\notag\\&+h(v)^{2}
\Bigg(c_{2} s_{0}(4B\eta -v^{2}\eta^{2}+2Bk_{0}\mu)+B^{2}\Bigg(-c_{1}^{2}s_{1}+k_{0}^{2}s_{1}+2c_{1}c_{2}s_{0}\Bigg)h^{\prime}(v)\Bigg),\notag\\
b_{0}(v)&=\frac{1}{4Bc_{2}h(v)}\Bigg(2c_{2}s_{0}\mu-2c_{2}s_{0}v\eta h(v)+Bc12s_{1}h^{2}(v)-Bk_{0}2s1h^{2}(v)\notag\\
&-2Bc_{1}c_{2}s_{0}h^{2}(v)-2Bc_{2}s_{0}h^{\prime}(v)\Bigg),\ \ \ c_{0}=\frac{c_{1}^{2}-k_{0}^{2}}{4c_{2}}.
\end{align}

The compatibly equation between \ $b_{0}^{\prime}(v)$ \ and \ $b_{0}(v)$
gives rise to

\begin{equation}
\begin{array}{c}
s_{0} \Bigg(\mu^{2}+\Bigg(-v^{2}\eta^{2}+2B(\eta+k_{0}\mu)\Bigg)h^{2}(v)+B^{2}k_{0}^{2}h^{4}(v)-4B\mu h^{\prime}(v)\\
+3B^{2}h^{\prime2}(v)-2B^{2}h(v)h^{\prime\prime}(v)\Bigg)=0.
\end{array}
\end{equation}

The equation has the first integral

\begin{equation}
h^{\prime}(v)=\frac{1}{B}(\eta+v\eta h(v)+Bk_{0}h^{2}(v)),
\end{equation}

that integrates to
\begin{equation}
\begin{array}{c}
h(v)=-\frac{P(v)}{Q(v)},\;P(v)=\mu(\sqrt{B}\sqrt{2}n\sqrt{\eta}B_{0}H_{n-1}\Bigg(\frac{v\sqrt{\eta}}{\sqrt{2B}}\Bigg)-nv\eta\:_{1}F_{1}(1-\frac{n}{2},\frac{3}{2},\frac{v^{2}\eta}{2B})),\\
Q(v)=n\,\eta B(B_{0}H_{n}(\frac{v\sqrt{\eta}}{\sqrt{2B}})-nv\eta\:_{1}F_{1}(-\frac{n}{2},\frac{1}{2},\frac{v^{2}\eta}{2B})),\ \ \ k_{0}=\frac{n\,\eta}{\mu},
\end{array}
\end{equation}

where $H_{n}(x)$ and $\ \ \   _{1}F_{1}(a,b,x)$ are the Hermite polynomial
and hyper-geometric function respectively.

The compatibly equation between $a_{1}^{\prime}(v)$ and $a_{1}(v)$
holds identically. It remains to find $a_{0}(v),$ by using (39),
we get

\begin{equation}
a_{:0}(v)=e^{\frac{v^{2}\eta}{2B}}B_{1},\ \ \ a_{1}(v)=\frac{2c_{2}a_{0}(v)}{c_{1}-k_{0}}.
\end{equation}

The solution of the auxiliary equation (38) gives rise to

\begin{equation}
g(v,t)=\frac{(n\,\eta-e^{n\,\eta(t+\int h(v)dv)}n\,\eta-c_{1}\mu-c_{1}e^{n\,\eta(t+\int h(v)dv)})\mu)}{2c_{2}\mu(1+e^{(n\,\eta(t+\int h(v)dv)})}.
\end{equation}

By substituting from (43) and (44) into the first equation in (37)
we have

\begin{equation}
\begin{array}{c}
f(t,v)=\frac{P_{1}}{Q_{1}},\;P_{1}=e^{\frac{v^{2}\eta}{2B}}(c_{1}-k_{0})((-1+e^{n\,\eta(t+\int h(v)dv)})ns_{1}\eta+\\
(1+e^{n\,\eta(t+\int h(v)dv)})(c_{1}s_{1}-2c_{2}s_{0})\mu),\;Q_{1}=2B_{1}c_{2}n((-1+e^{n\,\eta(t+\int h(v)dv)})n\eta\\
+(1+e^{n\,\eta(t+\int h(v)dv)})k_{0}\mu).
\end{array}
\end{equation}

Now we evaluate $\int h(v)dv$, to this end we use (40) and assume
that

\begin{equation}
P(v)=(Q(v)+S(v))^{\prime},
\end{equation}

where $S(v)$ is to be determined. Calculations show that

\begin{eqnarray}
S(v) & = & \frac{(n\eta+\mu)B}{(n\sqrt{\eta}}(n\sqrt{\eta})(-B_{0}H_{n}(\frac{v\sqrt{\eta}}{\sqrt{2B}})+(1-_{P}F_{Q}(\{-\frac{n}{2}\},\{\frac{1}{2}\},\frac{v^{2}\eta}{2B})),
\end{eqnarray}

where $_{P}F_{Q}$ $({a},{b},x)$ is the generalized hyper geometric
function {[}19{]}. Finally we get

\begin{equation}
\int h(v)dv=-Log(\mid Q(v)+S(v)\mid),
\end{equation}

where $Q(v)$ and $S(v)$ are give by (45) and (47) respectively.

\subsection{Self Similar Solution}

Now we consider the similarity transformations $z=v\,\omega(t),\ \ \ t:=t,$ \ and
$f(t,v)=\bar{f}(t,z)$ so the equation (2) is

\begin{equation}
\frac{\partial}{\partial t}\bar{f}(z,t)=\eta\frac{\partial}{\partial z}(z\bar{f}(z,t))+B\omega(t)^{2}\frac{\partial^{2}}{\partial z^{2}}\bar{f}(z,t),
\end{equation}

and as in subsection 4.1, (48) is transformed to

\begin{equation}
\begin{array}{c}
\tilde{f}_{z}(z,t)-F(z,t)\,\tilde{f}(z,t)=0,\:\tilde{f}_{t}(z,t)-G(z,t)\,\tilde{f}(z,t)=0,\\
G(z,t)-\eta(1+zF(z,t))-B\omega(t)^{2}(F_{z}(z,t)+F(z,t)^{2})=0.
\end{array}
\end{equation}

We assume that the solutions, by using (50), have the form

\begin{equation}
\begin{array}{c}
\tilde{f}(z,t)=\frac{s_{1}g(z,t)+s_{0}}{a_{1}(z,t)g(z,t)+a_{0}(z,t)},\:F(z,t)=\frac{b_{1}(z,t)g(z,t)+b_{0}(z,t)}{s_{1}g(z,,t)+s_{0}},\:\\
G(z,t)=\frac{d_{1}(z,t)g(z,t)+d_{0}(z,t)}{s_{1}g(z,,t)+s_{0}},
\end{array}
\end{equation}

and

\begin{equation}
\begin{array}{c}
g_{t}(z,t)=\mu(t)\,(c_{2}g(z,t)^{2}+c_{1}g(z,t)+c_{0}),\\
g_{z}(z,t)=h(z)\,(c_{2}g(z,t)^{2}+c_{1}g(z,t)+c_{0}).
\end{array}
\end{equation}
By inserting (51) and (52) into (50), we have following equations
\begin{equation}
\begin{array}{c}
a_{0z}(z,t)=\frac{1}{s_{0}\omega(t)}(-a_{0}(z,t)b(z,t)-c_{0}s_{0}a_{1}(z,t)h(z)\omega(t)\\
+c_{0}s_{1}a_{0}(z,t)\,h(z)\omega(t),\\
a_{1z}(z,t)=\frac{1}{s_{1}\omega(t)}(-a_{1}(z,t)b_{1}(z,t)-c_{2}s_{0}a_{1}(z,t_{0}\,h(z)\omega(t)\\
+c_{2}s_{1}a_{0}(z,t)\,h(z)\omega(t)),\\
b_{1}(z,t)=\frac{1}{s_{0}}(s_{1}b_{0}(z,t)-c_{0}s_{1}^{2}h(z)\omega(t)+c_{1}s_{1}s_{0}h(z)\omega(t)\\
-c_{2}s_{0}^{2}h(z)\omega(t)),
a_{0t}(z,t)=\frac{1}{s_{0}}(-a_{0}(z,t)d_{0}(z,t)-c_{0}s_{0}a_{1}(z,t)\mu(t)\\
+c_{0}s_{1}a_{0}(z,t)\mu(t)),\\
a_{1t}(z,t)=\frac{1}{s_{1}}(-a_{1}(z,t)d_{1}(z,t)-c_{2}s_{0}a_{1}(z,t)\mu(t)+c_{2}s_{1}a_{0}(z,t)\mu(t))\\
b_{0z}(z,t)=\frac{1}{4B^{2}h(z)^{2}\omega(t)^{3}}(-s_{0}\mu(t)^{2}+\\
2B\mu(t)\omega(t)^{2}((-2c_{0}s_{1}+c_{1}s_{0})h(z)^{2}+s_{0}h^{\prime}(z))\\
+\omega(t)^{2}(-B^{2}(c_{1}^{2}-4c2co)s_{0}h(z)^{4}\omega(t)^{2}-B^{2}s_{0}\omega(t)^{2}h^{\prime2}(z)\\
+h(z)^{2}(s_{0}\eta(-4B+z^{2}\eta)+4Bd_{0}(z,t)
\end{array}
\end{equation}
\begin{equation}
\begin{array}{c}
+2B^{2}(2c_{0}s_{1}-c_{1}s_{0})\omega(t)^{2}h^{\prime}(z))),\\
b_{0}(z,t)=\frac{1}{2Bh(z))\omega(t)}(s_{0}\mu(t)-s_{0}z\eta h(z)\omega(t)+2Bc_{0}s_{1}h(z)^{2}\omega(t)^{2}\\
-Bc_{1}s_{0}h(z)^{2}2\omega(t)^{2}-Bs_{0}\omega(t)^{2}h^{\prime}(z)).\notag
\end{array}
\end{equation}

Now the compatibly equation $b_{0z}(z,t)-(b_{0}(z,t))_{z}=0$ , solves
to

\begin{equation}
\begin{array}{c}
d_{0}(z,t)=\frac{1}{4Bh(z)^{2}\omega(t)^{2}}(s_{0}\omega(t)^{2}+2B\mu(t)\omega(t)^{2}((2c_{0}s_{1}-c_{1}s_{0})h(z)^{2}\\
-2s_{0}h^{\prime}(z))+s_{0}\omega(t)^{2}(B^{2}(c_{1}^{2}-4c_{2}c_{0})h(z)^{4}\omega(t)^{2}\\
-\eta h(z)^{2}(-4A+z2\eta+2B\omega(t))\\
+3B^{2}\omega(t)^{2}h^{\prime2}(z)-2B^{2}h(z)\omega(t)^{2}h^{\prime\prime}(z))).
\end{array}
\end{equation}

The compatibly equations \ $(a_{iz}(z,t))_{t}-(a_{it}(z,t))_{z}=0,$
give rise to

\begin{equation}
\begin{array}{c}
(-B^{2}(c_{1}^{2}-4c_{2}c_{0})h(z)^{4}\omega(t)^{5}h^{\prime}(z)+\omega(t)h(z)(\mu(t)^{2\mu(t)}\\
-4B\mu(t)\omega(t)^{2}h^{\prime}(z)+3B^{2}\omega(t)^{4}h^{\prime2}(z)+z\eta h(z)^{3}\omega(t)(\eta\omega(t)^{2}+\omega^{\prime}(t))\\
+2Bh(z)\omega(t)^{3}(\mu(t)-2B\omega(t)^{2}h^{\prime}(z))h^{\prime\prime}(z)\\
+h[z]2(\omega(t)\mu^{\prime}(t)-2\mu(t)\omega(t))+B^{2}\omega(t)^{5}h^{(3)}(z)))=0.
\end{array}
\end{equation}

We find that (54) holds when

\begin{equation}
\begin{array}{c}
\mu(t)=A_{0}\omega^{2}(t),\;\omega^{\prime}(t)=-\eta\omega(t)^{2},\ \ \ c_{0}=(c_{1}^{2}-p_{3}^{2})/(4c2),\\
h^{\prime}(z)=\frac{A_{0}}{A}+p_{1}h(z)+p_{2}h(z)^{2}.
\end{array}
\end{equation}

It remains to evaluate $a_{j}(z,t),\ \ \ j=0,1,...$, \ where detailed calculations
yield

\begin{equation}
\begin{array}{c}
a_{1}(z,t)=\frac{B_{3}\sqrt{B_{0}+t\eta}}{\sqrt{B_{0}}}e^{(\frac{1}{4}(2p_{1}z-4t\eta+(z^{2}\eta B_{0}+t\eta))/B-(Bp_{1}^{2}t)/(B_{0}^{2}+B_{0}t\eta)+(4A_{0}p_{2}t)/(B_{0}^{2}+Bot\eta)-4c_{1}\int h(z)dz))},\\
\\
a_{0}(z,t)=e^{\int_{0}^{t}\frac{1}{4}((\eta(-4+\frac{z^{2}\eta}{B})+2\eta\omega(t_{1})-Bp_{1}^{2}\omega(t_{1})^{2})dt_{1}}e^{(\frac{p_{1}z}{2}+\frac{B_{0}z^{2}\eta}{4B}}B_{1},\ \ \ c_{0}=0,\ \ \ c_{2}=0.
\end{array}
\end{equation}

The solution of the auxiliary equation is

\begin{equation}
g(z,t)=B_{2}e^{c_{1}(\int_{0}^{t}\mu(t_{1})dt_{1}+\int h(z)dz)},
\end{equation}

together with

\begin{equation}
h(z)=-\frac{-p1\sqrt{B}-rtanh(\frac{r(z+A_{1})}{2\sqrt{B}})}{2c_{1}\sqrt{B}}, \ \ \ r=\sqrt{Bp_{1}^{2}+4A_{0}c_{1}}.
\end{equation}

Finally we get the solution of (50) which is

\begin{equation}
\begin{array}{c}
f(z,t)=\frac{P_{2}}{Q_{2}},\ \ P_{1}=\sqrt{B_{0}}e^{-\frac{p_{1}z}{2}-\frac{z^{2}\eta(Bo+t\eta)}{4B}+t(\eta+\frac{Bp_{1}^{2}}{4B_{0}^{2}+4B_{0}t\eta})} \\
(B_{2}e^{\frac{1}{2}(\frac{(2A_{0}(-1+t)c_{1})}{(B_{0}+\eta)(B_{0}+t\eta)}-p_{1}z-2Log(cosh((\frac{\sqrt{Bp_{1}^{2}+4A_{0}c_{1}}(z+A_{1})}{2\sqrt{B}}))}s_{1}+s_{0}),\\
Q_{1}=(B_{1}+B_{2}B_{3}e^{\frac{A_{0}(B_{0}(c_{1}(-1+t)-c_{1}t)-c_{1}t\eta)}{B_{0}(B_{0}+\eta)(Bo+t\eta)}})\sqrt{Bo+t\eta},\  \ \ z=v,\ \ \ \omega(t)=\frac{v}{\sqrt{Bo+t\eta}}.
\end{array}
\end{equation}

\section{Solutions of Fractional FPE}

The study of fractional evolution equations occupies a remarkable
area in the literature {[}29-32{]}.  Here the fractional
time derivative FPE is reduced to the classical one's by using the
similarity transformations $f(t,v)=\tilde{f}(\tau, v)$ \ and \  $\tau$
is given by


$(a_{1})$ in the case CFD $\tau=\frac{\Gamma(2-\alpha)}{\alpha}(T_{0}^{\alpha}-(T_{0}-t)^{\alpha})$.

$(a_{2})$ in the case of CFFD $\tau=\frac{(2-\alpha)(1-\alpha)}{2\alpha}Log(\frac{e^{\frac{\alpha}{1-\alpha}(T_{0}-t)}-1}{e^{\frac{\alpha}{1-\alpha}T_{0}}-1}).$

$(a_{3})$ in the case of  GFD $\tau=\frac{(\lambda+2)}{2}\int_{0}^{t}\frac{1}{\gamma(\frac{1}{\beta},\lambda(T_{0}-t_{1})^{\beta})}dt_{1}.$

The equation (2) is rewritten as follows:

\begin{equation}
\frac{\partial}{\partial\tau}\widetilde{f}(\tau,v)=\eta\frac{\partial}{\partial v}(v\tilde{f}(\tau,v))+B\frac{\partial^{2}}{\partial v^{2}}\tilde{f}(\tau, v).
\end{equation}

We find two classes of solutions, first by using the similarity transformations,
where (61) is used and second by considering (2) as a non autonomous
equation

\begin{equation}
p(t)\frac{\partial}{\partial t}f(t,v)=\eta\frac{\partial}{\partial v}(vf(t,v))+B\frac{\partial^{2}}{\partial v^{2}}f(t,v),
\end{equation}
where $p(t)$ takes one of the forms:


$(b_{1})$ in the case CFD $p(t)=\frac{1}{\Gamma(2-\alpha)}(T_{0}-t)^{1-\alpha}$.

$(b_{2})$ in the case of CFFD $p(t)=\frac{2}{(2-\alpha)}(1-e^{-\frac{\alpha}{1-\alpha}(T_{0}-t)}).$

$(b_{3})$ in the case of GFD $p(t)=\frac{2}{\lambda+2}\gamma(\frac{1}{\beta},\lambda(T_{0}-t)^{\beta}).$

Here, we consider (61) and by using $(a_{i}),i=1,2,3,4$, the numerical
results of the solution, which given by (35) and (36), are evaluated.

We confine ourselves to present some numerical results with relevance
to subsection 4.1. Attention is focused to consider the cases of fractional
and CFD.

The results of the solution, given in subsection 4.1 are displayed
against $v$ and $t$, and they are shown in Figures 1 (i)-(iii) for
different values of  $\beta$  and the friction coefficient $\eta$.

\begin{figure}[htbp]
	\includegraphics[width=1\textwidth]{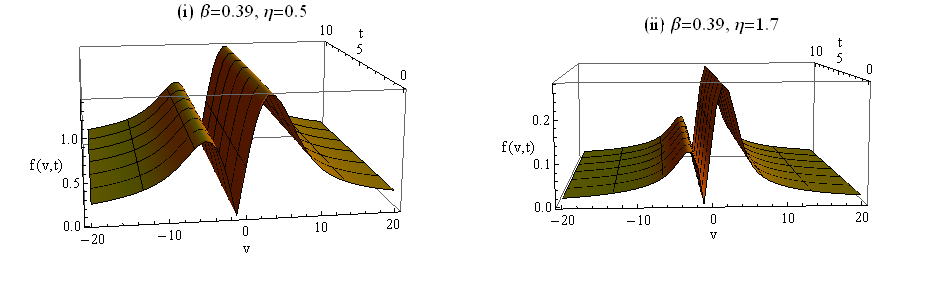}
	\begin{center}
		\includegraphics[width=0.5\textwidth]{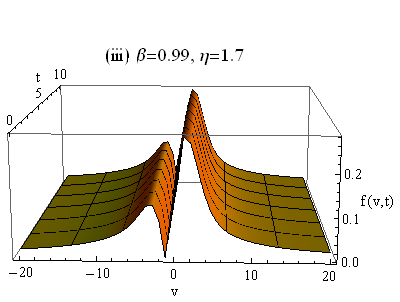}
	\end{center}
\caption{\small{ The behavior of  distribution function, when (i) $B_{0}=1.5, \ c_{0}=2, \ n=10, \ \eta=1.7, \ A_{0}=1.3, \ A_{1}=2.3, \ B=5, \ m=2.5, \ \tau=t^{\beta}, \ \beta=0.39, \ A_{2}=3, \ \mu=-0.5, \ B_{1}=1.9, \ B_{2}=0.7$.
		(ii) the same caption as in (i) but $\eta=1.7.$ \ (iii) the same caption
		as in (i) but  \ $\beta=0.99.$ }}	
\end{figure} 
Figures 1  show that the distribution function is mixed-Gaussian's
and that the friction coefficient plays a significant role in lowering
the magnitude of the distribution density function. While the
effect of vanning the fractional order plays a in lowering the tails. 
\begin{figure}[htbp]
	\includegraphics[width=1\textwidth]{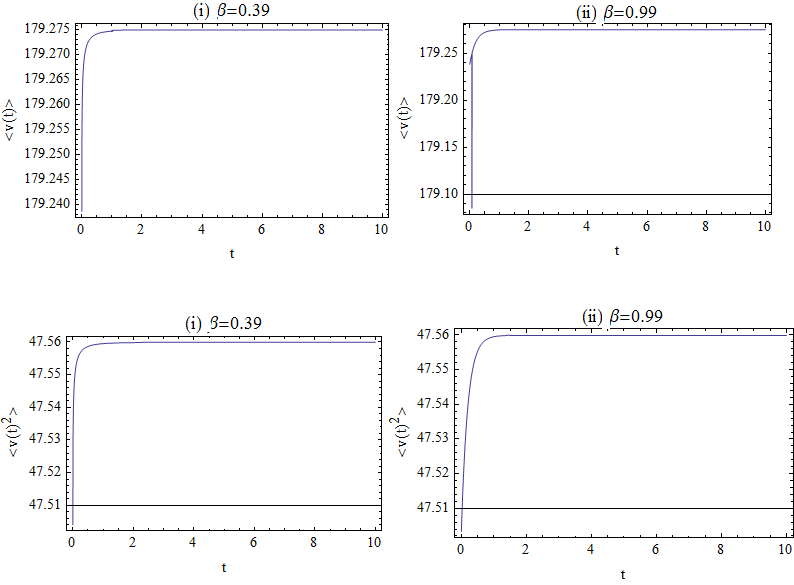}
		\caption{\small{ The mean and  mean square of the velocity,
				when $B_{0}=1.5, \ c_{0}=2, \ n=10, \ \eta=0.5, \ A_{0}=1.3, \ A_{1}=2.3, \ B=5, \ m:=2.5, \ \tau=t^{\beta}, \ \beta=0.99, \ A_{2}=3, \ \mu=-0.5, \ B_{1}=1.9, \ B_{2}=0.7.$
		}}
\end{figure} 



The order of the fractional time derivative has no remarkable effect
on the mean and the mean square.

For the solution of fractional FPE in the caputo sense, the results
in subsection 5.1 are displayed against $v$ and $t$, and they are shown
in Figures3 (i)-(iii), for different values of $\alpha$ and the friction
coefficient  \ $\eta$.
\begin{figure}[htbp]
	\includegraphics[width=1\textwidth]{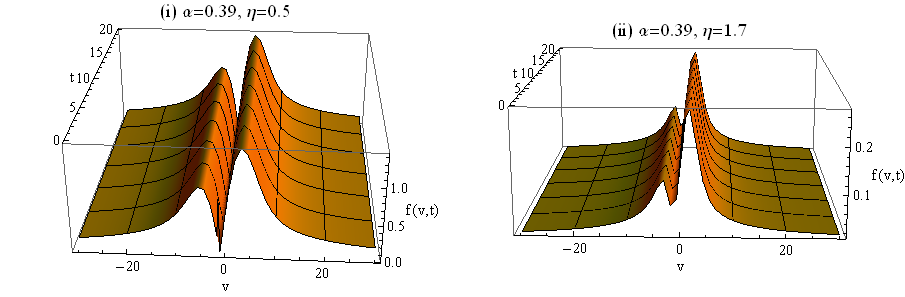}
		\begin{center}
			\includegraphics[width=0.5\textwidth]{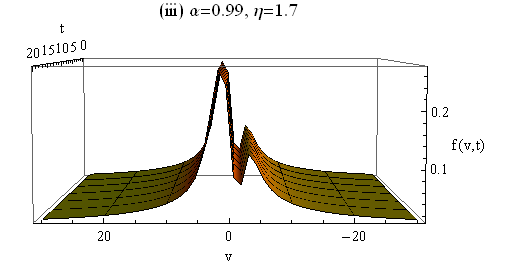}
			\caption{\small{ The behavior of  distribution function, when (i) $B_{0}=1.5, \ c_{0}=2, \ n=10, \ \eta=0.5, \ A_{0}:=1.3, \ A_{1}=2.3, \ B=5, \ m=2.5, \ \tau=\frac{\Gamma(2-\alpha)}{\alpha}(T_{0}^{\alpha}-(T_{0}-t)^{\alpha}), \ \ \alpha=0.39, \ T_{0}=20, \ A_{2}:=3, \ \mu=-0.5, \ B_{1}:=1.9, \ B_{2}:=0.7.$
					(ii) the same caption as in (i) but $\eta=1.7.$ \ (iiii) the same caption
					as in (i) but $\alpha=0.99.$}}
		\end{center}	
\end{figure} 

%


Figures 3 show mixed Gaussian's. The friction coefficient plays a
dominant role in lowering the magnitude of the distribution function.
While when $\alpha=0.99$, permutation the Gaussian's occurs.

Figures 4 show  the mean and mean square of the velocity
are displayed against t, by varying the fractional order.
\begin{figure}[htbp]
	\includegraphics[width=1\textwidth]{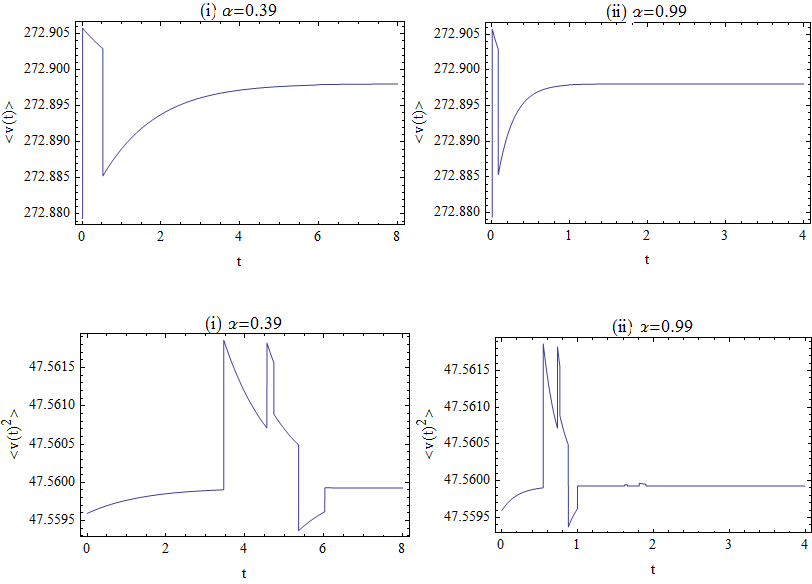}
	\caption{\small{ The mean and mean square for the same caption as in
			Figs.3.
	}}
	
\end{figure} 


Figures 4 (i) and (ii), show no remarkable variation in the mean and
mean square when varying the fractional order.

\section{Conclusions}

An approach for finding solutions of linear PDEs with variable coefficients
is presented. It is established by transforming the PDE to a system
of first order PDEs and the extended unified method is implemented
a class of solutions of  fractional Fokker Planck equations
are obtained. The solutions show that the distribution function is
mixed-Gaussian's. Further the friction coefficient plays the role
of lowering the magnitude of the distribution function. On the other
hand, varying the order of the fractional time derivative has the
effect of permuting the Gaussian's the distribution function.

\end{document}